\documentclass[final,nomarks]{dmtcs-episciences}
\usepackage{amsmath}
\usepackage{amssymb}
\usepackage{amsthm}
\usepackage{hyperref}
\usepackage[capitalise,noabbrev]{cleveref}
\usepackage{booktabs}
\usepackage{multirow}
\usepackage{color}
\usepackage{enumitem}

\usepackage[round]{natbib}
\usepackage{doi}

\newtheorem{thm}{Theorem}[section]
\newtheorem{prop}[thm]{Proposition}

\newtheorem{coro}[thm]{Corollary}
\theoremstyle{remark}
\newtheorem{rmk}[thm]{Remark}

\theoremstyle{definition}
\newtheorem{defn}[thm]{Definition}
\newtheorem{constr}[thm]{Construction}

\crefname{thm}{Theorem}{Theorems}
\crefname{prop}{Proposition}{Propositions}
\crefname{lem}{Lemma}{Lemmas}
\crefname{coro}{Corollary}{Corollaries}
\crefname{rmk}{Remark}{Remarks}
\crefname{example}{Example}{Examples}
\crefname{conj}{Conjecture}{Conjectures}
\crefname{defn}{Definition}{Definitions}
\crefname{constr}{Construction}{Constructions}

\newcommand{\tdef}[1]{\textcolor{blue}{\emph{#1}}}

\newcommand{\insertfig}[2]{\includegraphics[page=#1, width=#2\textwidth]{fig-ipe.pdf}}

\newcommand{\dinv}{\mathtt{dinv}}
\newcommand{\bounce}{\mathtt{bounce}}

\newcommand{\luka}{Łukasiewicz}
\newcommand{\tC}{\widetilde{C}}
\newcommand{\lukaset}{\mathcal{L}}
\newcommand{\trees}{\mathcal{T}}

\newcommand{\profile}{\mathtt{Pr}}
\newcommand{\profmset}{\mathtt{M}}
\newcommand{\area}{\mathtt{area}}
\newcommand{\areavec}{\mathtt{Area}}
\newcommand{\depth}{\mathtt{depth}}
\newcommand{\depthvec}{\mathtt{Depth}}
\newcommand{\lthorn}{\mathtt{lthorn}}
\newcommand{\rthorn}{\mathtt{rthorn}}

\newcommand{\lukatotree}{\tau}
\newcommand{\treetoluka}{\lambda}
\newcommand{\mirror}{\operatorname{mir}}
\newcommand{\swaplode}{\operatorname{swap}}

\author{Wenjie Fang}

\affiliation{Univ Gustave Eiffel, CNRS, LIGM, F-77454 Marne-la-Vallée, France}


\title{On the area-depth symmetry on \luka{} paths}

\begin{document}

\publicationdata{vol. 28:3}{2026}{2}{10.46298/dmtcs.17405}{2026-01-27; 2026-01-27; 2026-05-29}{2026-06-22}

\maketitle

\begin{abstract}
  \vspace*{2ex} % problem with dmtcs template
  In an effort to further understanding $q,t$-Catalan statistics, a new statistic on Dyck paths called $\depth$ was proposed in Pappe, Paul and Schilling (2022) and was shown to be jointly equi-distributed with the well-known $\area$ statistics. In a recent preprint, Qu and Zhang (2025) generalized $\depth$ to so-called ``$\vec{k}$-Dyck paths''. They showed that $\area$ and $\depth$ are also jointly equi-distributed over such paths with a fixed multiset of up-steps and a given first up-step, and they conjectured that the same holds when also fixing the last up-step. In this short note, we settle this conjecture on the more general context of \luka{} paths by interpreting $\area$ and $\depth$ under the classical bijection between \luka{} paths and plane trees, through which the symmetry is transparent.
\end{abstract}

\section{Introduction} \label{sec:intro}

The study of $q,t$-Catalan polynomials $C(q, t)$ started from two combinatorial formulas of the bi-graded Hilbert series of the subspace of alternants in the diagonal coinvariant space with two sets of variables, both describing Dyck paths of a given size with two jointly equi-distributed statistics: with $\area$ and $\bounce$ by \citet{area-bounce}, and with $\dinv$ and $\area$ by Haiman. The two pairs of statistics are related by the famous zeta map \citep*{zeta-map} (see also \citep*{cataland} for a simpler definition). For more details, readers are referred to \cite{haglund-book}.

From the algebraic definition, it is trivial that the $q,t$-Catalan polynomials are symmetric in $q$ and $t$. However, on the combinatorial side, it is totally mysterious, with no bijective explanation known. In an effort to better understand combinatorially the related statistics, \citet*{dyck-depth} proposed another statistic $\depth$ on Dyck paths and found that it is also jointly equi-distributed with $\area$ by exhibiting an involution exchanging $\area$ and $\depth$. The $q,t$-polynomial $\tC(q, t)$ defined with $\area$ and $\depth$, which is different from the classical $C(q, t)$, is thus also $q,t$-symmetric.

There are many generalizations of the $q,t$-Catalan polynomials, usually defined on generalizations of Dyck paths, sometimes with an algebraic motivation, but not always. For instance, \citet{kvec-stats} considered in the so-called \emph{$\vec{k}$-Dyck paths}, on which they introduced the three classical statistics $\dinv$, $\area$ and $\bounce$, although without direct algebraic supporting background. They then defined a generalized $q,t$-Catalan polynomial $C_\lambda(q, t)$ on $\vec{k}$-Dyck paths using these generalized statistics and studied their properties. Unfortunately, $C_\lambda(q, t)$ is not $q,t$-symmetric in general.

In search for generalizations of the $q,t$-Catalan polynomials defined on $\vec{k}$-Dyck paths with better $q,t$-symmetric properties, \citet*{kvec-depth} generalized the $\depth$ statistics from \citet*{dyck-depth} to $\vec{k}$-Dyck paths. They defined a generalization of the $\area$-$\depth$ $q, t$-polynomial $\tC(q, t)$ to $\vec{k}$-Dyck paths with a given multiset $M$ of up-steps, denoted here by $\tC_M(q, t)$. They also defined variants of $\tC_M(q, t)$ with extra conditions on the first and/or the last up-step. They then showed bijectively that the variant of $\tC_M(q, t)$ fixing the first up-step is $q, t$-symmetric, and they conjectured that it remains $q, t$-symmetric even when also fixing the last up-step. The main result of this short note is a proof of this conjecture for \luka{} paths, which generalize $\vec{k}$-Dyck paths by allowing horizontal steps. Precise definitions of notations are postponed to later sections.

\begin{thm}[{\cite[Conjecture~1.4]{kvec-depth}}] \label{thm:qt-sym-lodestar}
  We define $\tC_{a,M,b}(q, t)$ as a sum over \luka{} paths $P$ with the first (resp. last) up-step of degree $a$ (resp. $b$) and $M$ the multiset of the degrees of other up-steps:
  \[
    \tC_{a,M,b}(q, t) = \sum_P q^{\area(P)} t^{\depth(P)}.
  \]
  Then $\tC_{a,M,b}(q, t)$ is symmetric in $q, t$, \emph{i.e.}, $\tC_{a, M, b}(q, t) = \tC_{a, M, b}(t, q)$.
\end{thm}

Our proof relies essentially on the observation that both statistics $\area$ and $\depth$ on \luka{} paths can be interpreted naturally over the classical bijection between \luka{} paths and plane trees.

The rest of this article is organized as follows. We lay down the definitions of \luka{} paths, the statistics $\area$ and $\depth$ and the related polynomials in \Cref{sec:luka}. Then we revisit the classical bijection between \luka{} paths and plane trees in \Cref{sec:trees}, along with how $\area$ and $\depth$ are transferred to plane trees. Finally, in \Cref{sec:sym} we show how to obtain $q,t$-symmetry results using simple bijections in a transparent way.

\paragraph{Acknowledgment} We would like to thank the anonymous reviewers for their helpful comments. This work is not supported by any funding with a predefined goal, but by the publicly funded laboratory LIGM of Université Gustave Eiffel.

\section{\luka{} paths and its statistics} \label{sec:luka}

\begin{figure}
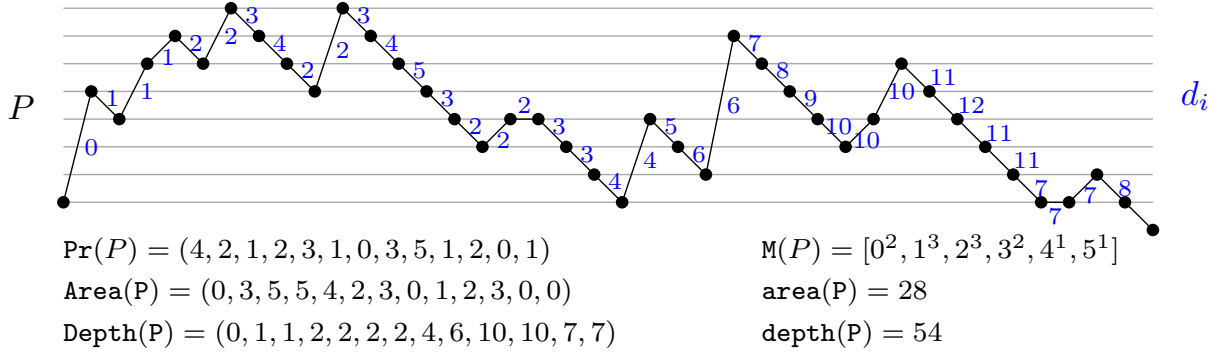

  \centering
  \insertfig{1}{1}
  \caption{Example of a \luka{} path, with its statistics $\area$ and $\depth$.}
  \label{fig:luka}
\end{figure}

\tdef{\luka{} paths} are lattices paths in $\mathbb{Z}^2$ with steps of the form $(1, k)$ with $k \geq -1$ that go under the $x$-axis only at the last step, which is always $(1, -1)$. We denote by $\lukaset$ the set of \luka{} paths. \Cref{fig:luka} shows an example of \luka{} paths. For the ease of notation, we denote by $D$ the only down-step $(1, -1)$, and $U_k$ the up-step $(1, k)$ for $k \geq 0$, and we say that $U_k$ is of \tdef{degree} $k$. The \tdef{degree profile} (or simply \tdef{profile}) of a \luka{} path $P$, denoted by $\profile(P)$, is the sequence of degrees of up-steps in $P$ from left to right. In this case, the \tdef{profile multiset} of $P$, denoted by $\profmset(P)$, is the multiset of entries in $\profile(P)$. For a finite multiset $M$ of natural numbers, we denote by $M_k$ the multiplicity of $k$ in $M$, and we write $M = [0^{M_0}, 1^{M_1}, 2^{M_2}, \dots]$ while only listing elements in $M$. We note that $\vec{k}$-Dyck paths, the family of paths studied by \citet*{kvec-depth}, first introduced by \citet*{kvec-sweep}, are simply \luka{} paths with $\profile(P) = \vec{k}$ without any $0$ and with the last down-step removed.

We now introduce the statistics we consider on \luka{} paths, which is a simple generalization of those from \cite{kvec-depth} for $\vec{k}$-Dyck paths.

\begin{defn} \label{defn:luka-area}
Given a \luka{} path $P$, its \tdef{area vector} $\areavec(P)$ is the vector whose $i$-th entry is the $y$-coordinate of the starting point of the $i$-th up-step of $P$. The \tdef{area} of $P$, denoted by $\area(P)$, is the sum of all entries in $\areavec(P)$. See \Cref{fig:luka} for an example.  
\end{defn}

\begin{defn} \label{defn:luka-depth}
In a \luka{} path $P$, we associate each down-step $P_j = D$ except the last one with an up-step $P_i = U_k$ with $k > 0$ such that, when drawing a horizontal ray from the midpoint of $P_j$ to the left, the first non-horizontal up-step it hits in the middle is $P_i$. In this case, we say that $P_j$ is \tdef{matched} with $P_i$. An up-step $P_j = U_k$ has exactly $k$ matching down-steps $P_{j_1}, \dots, P_{j_k}$ with $j_1 < \dots < j_k$, and we say that $P_{j_i}$ is the $i$-th matching down-step of $P_j$. We then associate to each step $P_i$ a value $d_i$ as follows:
\begin{itemize}
\item The first step $P_1$, which is always an up-step, has value $d_1 = 0$;
\item If $P_i = D$, let $P_j$ be its matching up-step, and suppose that $P_i$ is the $\ell$-th matching down-step of $P_j$, then we take $d_i = d_j + \ell$;
\item If $P_i = U_k$ with $i > 0$, we take $d_i = d_{i-1}$.
\end{itemize}
The \tdef{depth vector} of $P$, denoted by $\depthvec(P)$, is the vector formed by the values $d_i$ for all up-steps $P_i$ from left to right. The \tdef{depth} of $P$, denoted by $\depth(P)$, is the sum of all entries in $\depthvec(P)$. See \Cref{fig:luka} for an example.
\end{defn}

\begin{rmk}
  We note that our definition of $\depth$ seems different from the one in \cite{kvec-depth}, which used two tableaux called ``Filling tableau'' and ``Ranking tableau'', whose columns represent up-steps in a $\vec{k}$-Dyck path $P$. In a word, the Filling tableau of $P$ can be seen as associating $k$ down-steps to each up-step $P_i = U_k$ of degree $k$ by scanning down-steps from left to right and associate each down-step with the closest up-step before it that still needs new down-steps. This is clearly the same as our notion of matching steps. The Ranking tableau pf $P$ is then obtained by labeling the first up-step by $0$, then the $i$-th down-step associated to an up-step with label $\ell$ is labeled by $\ell + i$, while a remaining up-step is also labeled $\ell$ if its previous step is labeled $\ell$. The depth of $P$ is the sum of labels of all up-steps. This also agrees with our definition here.
\end{rmk}

Now we define the $q,t$-polynomial associated to $\area$ and $\depth$. We denote the concatenation of two vectors $P, Q$ by $P \cdot Q$, and the union of two multisets $L, M$ by $L \uplus M$. We further give four notations of subsets of \luka{} paths:
\begin{itemize}
\item $\lukaset_M$: the set of \luka{} paths $P$ with $\profmset(P) = M$;
\item $\lukaset_{a, M}$: the set of \luka{} paths $P$ with $U_a$ as the first up-step, and $\profmset(P) = M \uplus \{a\}$;
\item $\lukaset_{M, b}$: the set of \luka{} paths $P$ with $U_b$ as the last up-step, and $\profmset(P) = M \uplus \{b\}$;
\item $\lukaset_{a, M, b}$: the set of \luka{} paths $P$ with $U_a$ (resp. $U_b$) as the first (resp. last) up-step, and $\profmset(P) = M \uplus \{a, b\}$, with $\{a, b\}$ understood as a multiset.
\end{itemize}
We then define
\begin{align}
  \begin{split}
    \label{eq:qt-poly-def}
    \tC_{M}(q, t) = \sum_{P \in \lukaset_{M}} q^{\area(P)} t^{\depth(P)}, &\quad \tC_{a, M}(q, t) = \sum_{P \in \lukaset_{a, M}} q^{\area(P)} t^{\depth(P)}, \\
    \tC_{M, b}(q, t) = \sum_{P \in \lukaset_{M, b}} q^{\area(P)} t^{\depth(P)}, &\quad \tC_{a, M, b}(q, t) = \sum_{P \in \lukaset_{a, M, b}} q^{\area(P)} t^{\depth(P)}.
  \end{split}
\end{align}

We may also define a more refined version of these polynomials. Given a vector $\mathbf{K}$, we define the following polynomial summing over \luka{} paths with profile $\mathbf{K}$.
\begin{equation}
  \label{eq:qt-refined-def}
  \tC_{\mathbf{K}}(q, t) = \sum_{P \in \lukaset, \profile(P) = \mathbf{K}} q^{\area(P)} t^{\depth(P)}.
\end{equation}
In the following, the profile vector will always be given in boldface to avoid confusion between the polynomials.

\section{Plane trees and $q,t$-statistics} \label{sec:trees}

A \tdef{plane tree} is rooted trees in which each \tdef{internal node}, \textit{i.e.}, a node with at least one child, comes with a left-to-right order of its children. The \tdef{degree} of an internal node is the number of its children. Nodes that are not internal are called \tdef{leaves}. We denote by $\trees$ the set of plane trees. There is a classical bijection between \luka{} paths and plane trees, see, \textit{e.g.}, \cite[Section~I.5.3]{flajolet}, that we describe below.

\begin{constr} \label{constr:luka-tree}
  Given a plane tree $T$, we perform a contour walk from left to right, starting from the root, and transcribe the nodes into steps upon their first encounter: $D$ if it is a leaf, and $U_k$ if it is an internal node of degree $k + 1$. The lattice path thus obtained, denoted by $\treetoluka(T)$, is a \luka{} path.

  Conversely, given a \luka{} path $P$, starting with an empty positioning called \emph{bud}, we construct a plane tree by reading $P$ from left to right: when reading an up-step $U_k$ (resp. down-step $D$), we replace the leftmost empty bud by an internal node with $k + 1$ empty buds below (resp. a leaf). The plane tree thus obtained is denoted by $\lukatotree(P)$. It is clear that $\treetoluka$ and $\lukatotree$ are bijections between $\lukaset$ and $\trees$, and are inverse of each other.
\end{constr}

\begin{figure}
  \centering
  \insertfig{2}{1}
  \caption{Example of the bijections $\treetoluka$ and $\lukatotree$.}
  \label{fig:luka-tree-bij}
\end{figure}

The statistics $\area$ and $\depth$ on \luka{} paths are transferred naturally to plane trees. We see in \Cref{constr:luka-tree} that up-steps of a \luka{} path $P$ from left to right are converted to internal nodes of $T = \lukatotree(P)$ in the contour walk order, with $U_k$ translated to a node of degree $k + 1$. In the following, we denote by $u_i$ the internal node in $T$ associated to the $i$-th up-step of $P$.

\begin{figure}
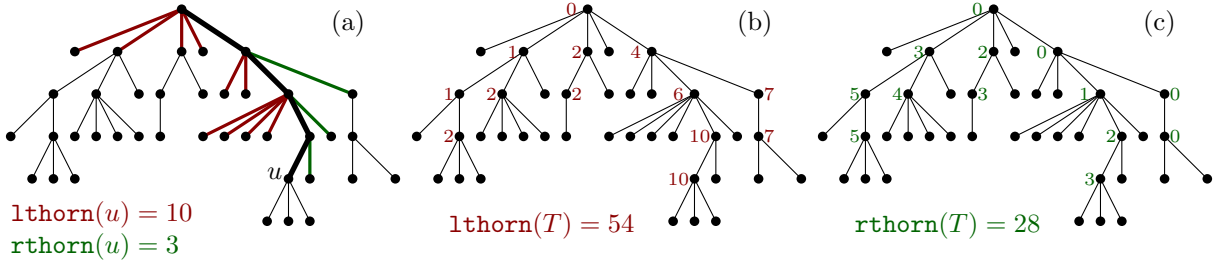

  \centering
  \insertfig{3}{1}
  \caption{(a) Illustration of $\lthorn$ and $\rthorn$ on an internal node. (b) Example of $\lthorn$. (c) Example of $\rthorn$.}
  \label{fig:tree-thorns}
\end{figure}

\begin{defn} \label{defn:tree-thorns}
  Given a plane tree $T$, for an internal node $u$, an edge $e$ is called a \tdef{left thorn} (resp. \tdef{right thorn}) of $u$ if $e$ is a descending edge of an ancestor $v \neq u$ to the left (resp. right) of the path from $v$ to $u$. We denote by $\lthorn(u)$ (resp. $\rthorn(u)$) the number of left thorns (resp. right thorns) of $u$, and by $\lthorn(T)$ (resp. $\rthorn(T)$) the sum of $\lthorn(u)$ (resp. $\rthorn(u)$) for all internal nodes $u$ of $T$. See \Cref{fig:tree-thorns} for an illustration.
\end{defn}

Comparing the statistics in \Cref{fig:luka,fig:tree-thorns}, it seems that $\area$ and $\depth$ on \luka{} paths are transferred to $\rthorn$ and $\lthorn$ on plane trees. This is indeed the case.

\begin{prop} \label{prop:area-rthorn}
  Let $P$ be a \luka{} path and $T = \lukatotree(P)$. Then $\area(P) = \rthorn(T)$.
\end{prop}
\begin{proof}
  We show instead the stronger result $\areavec(P)_i = \rthorn(u_i)$, with $u_i$ the internal node of $T$ associated to the $i$-th up-step of $P$ in \Cref{constr:luka-tree}. Our original claim comes from summing on both sides. We proceed by induction on the distance from the root to $u_i$. For distance $0$, we have $i = 1$ and $u_1$ is the root. It is clear that $\areavec(P)_1 = 0 = \rthorn(u_1)$. For non-zero distance, $u_i$ is a child of some internal node $u_j$, and by induction hypothesis we have $\areavec(P)_j = \rthorn(u_j)$. Suppose that $u_j$ has $k$ children, with $u_i$ the $\ell$-th one. By \Cref{defn:tree-thorns}, we have $\rthorn(u_i) = \rthorn(u_j) + k - \ell$. Now, in \Cref{constr:luka-tree}, it is clear that between the ending point of the $j$-th up-step of $P$ and the starting point of the $i$-th one is a concatenation of $\ell - 1$ \luka{} paths, each corresponding to the sub-trees of $T$ induced by the $\ell - 1$ children of $u_j$ to the left of $u_i$. As each \luka{} path introduces a variation of $-1$ in the $y$-coordinate, we have $\areavec(P)_i = \areavec(P)_j + (k - 1) - (\ell - 1)$. We thus have $\areavec(P)_i = \rthorn(u_i)$, by which we conclude the induction.
\end{proof}

\begin{prop} \label{prop:depth-lthorn}
  Let $P$ be a \luka{} path and $T = \lukatotree(P)$. Then $\depth(P) = \lthorn(T)$.
\end{prop}
\begin{proof}
  We show instead the stronger result $\depthvec(P)_i = \lthorn(u_i)$, with $u_i$ the internal node of $T$ associated to the $i$-th up-step of $P$ in \Cref{constr:luka-tree}. Our original claim comes from summing on both sides. We proceed by induction on the root distance of $u_i$. When $u_i$ is the root, we have $i = 1$, and it is clear that $\depthvec(P)_1 = 0 = \lthorn(u_1)$. In the case when $u_i$ has a parent $u_j$, by induction hypothesis we have $\depthvec(P)_j = \lthorn(u_j)$. Suppose that $u_j$ has $k$ children, with $u_i$ the $\ell$-th one. By \Cref{defn:tree-thorns} we have $\lthorn(u_i) = \lthorn(u_j) + \ell - 1$. Now, as in a \luka{} path, the only unmatched down step is the last one, it is clear that the $j$-th up-step of $P$ is matched with the down steps giving the rightmost leaf of sub-trees of $T$ induced by children of $u_j$, except for the last one. As the step preceding the $i$-th up-step of $P$ is the down step giving the rightmost leaf of the sub-tree induced by the $(\ell - 1)$-st child of $u_j$ (or the $j$-th up-step if $\ell = 1$), by \Cref{defn:luka-depth}, we have $\depthvec(P)_i = \depthvec(P)_j + \ell - 1$. We thus have $\depthvec(P)_i = \lthorn(u_i)$, by which we conclude the induction.
\end{proof}

\section{Symmetries  in $q,t$ through trees} \label{sec:sym}

Given a plane tree $T$, we denote by $\mirror(T)$ its image by left-right reflection. The following proposition is clear from \Cref{defn:tree-thorns}.

\begin{prop} \label{prop:thorns-sym}
  For a plane tree $T$, let $T' = \mirror(T)$, then $\lthorn(T) = \rthorn(T')$ and $\rthorn(T) = \lthorn(T')$.
\end{prop}
\begin{proof}
  It is clear that a left thorn $e$ of a node $u$ is a right thorn of $u$ in $\mirror(T)$.
\end{proof}

Now the $q,t$-symmetry of $\tC_{a, M}(q, t)$ is transparent.

\begin{thm}[{\cite[Theorem~1.1]{kvec-depth}}] \label{prop:qt-sym-plain}
  For any $a \in \mathbb{N}$ and multiset $M$ with elements in $\mathbb{N}$, we have $\tC_{a, M}(q, t) = \tC_{a, M}(t, q)$.
\end{thm}
\begin{proof}
  Let $P$ be a \luka{} path, and $T = \lukatotree(P)$. From \Cref{constr:luka-tree}, $P \in \lukaset_{a, M}$ if and only if $T$ has a root of degree $a + 1$, and the multiset of the degrees of non-root internal nodes is $M' = [1^{M_0}, 2^{M_1},\dots]$. Thus $P' = \treetoluka(T')$ with $T' = \mirror(T)$ is in $\lukaset_{a, M}$ if and only if $P \in \lukaset_{a, M}$, meaning that $\treetoluka \circ \mirror \circ \lukatotree$ is an involution on $\lukaset_{a, M}$. We then conclude by \Cref{eq:qt-poly-def,prop:area-rthorn,prop:depth-lthorn,prop:thorns-sym}.
\end{proof}

To show the $q,t$-symmetry of $\tC_{a, M, b}(q, t)$, we need an extra involution.

\begin{figure}
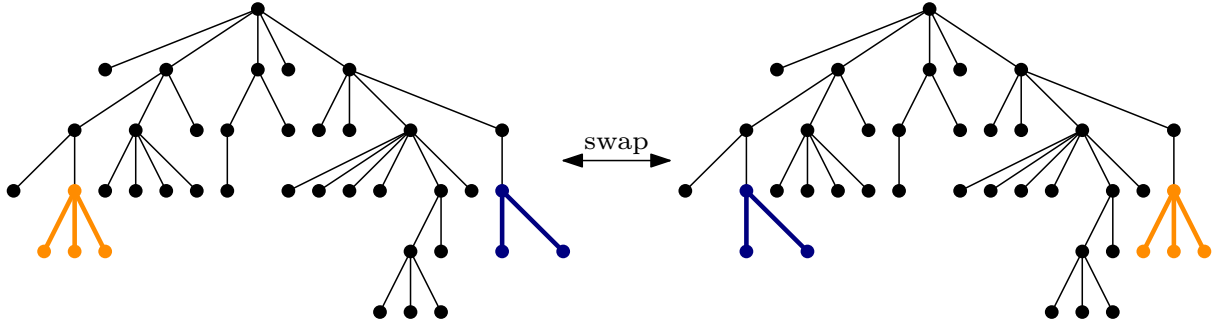

  \centering
  \insertfig{4}{1}
  \caption{Example of the lodestar swapping map.}
  \label{fig:lodestar}
\end{figure}

\begin{constr} \label{constr:swap-lodestar}
  Given a plane tree $T$, it \tdef{left lodestar} (resp. \tdef{right lodestar}) is the leftmost (resp. rightmost) internal node whose children are all leaves. Note that the left (resp. right) lodestar of $T$ is not necessarily on its leftmost (resp. rightmost) path, and it is possible that both lodestars are the same. We define the \tdef{lodestar swap} of $T$, denoted by $\swaplode(T)$, to be the tree obtained by swapping the left and right lodestars of $T$. See \Cref{fig:lodestar} for an example of the lodestar swapping map.
\end{constr}

\begin{prop} \label{prop:lodestar-sym}
  Let $T$ be a plane tree and $T' = \swaplode(T)$ its lodestar swap. Then $\lthorn(T) = \lthorn(T')$ and $\rthorn(T) = \rthorn(T')$.
\end{prop}
\begin{proof}
  By \Cref{defn:tree-thorns}, the descendant edges of lodestars are never thorns of any internal node, thus their swapping does not affect $\lthorn$ and $\rthorn$ of any node.
\end{proof}

Now we can show our main result \Cref{thm:qt-sym-lodestar}.

\begin{proof}[{of \Cref{thm:qt-sym-lodestar}}]
  Let $P$ be a \luka{} path, and $T = \lukatotree(P)$. It is clear from \Cref{constr:luka-tree} that the right lodestar of $T$ corresponds to the last up-step of $P$. From \Cref{constr:luka-tree}, $P \in \lukaset_{a, M, b}$ if and only if $T$ has a root of degree $a + 1$, the multiset of the degrees of non-root internal nodes is $M' = [1^{M_0}, 2^{M_1},\dots] \uplus \{b + 1\}$, and the right lodestar of $T$ is of degree $b + 1$. We take $P' = \treetoluka(T')$ with $T' = \swaplode(\mirror(T))$. It is clear that the left (resp. right) lodestar of $\mirror(T)$ is from the right (resp. left) lodestar of $T$. Hence, the right lodestar of $T'$ is from the right lodestar of $T$. We thus know that $P'$ is also in $\lukaset_{a, M, b}$ if and only if $P \in \lukaset_{a, M, b}$, meaning that $\treetoluka \circ \swaplode \circ \mirror \circ \lukatotree$ is an involution on $\lukaset_{a, M, b}$. We then conclude by \Cref{eq:qt-poly-def,prop:area-rthorn,prop:depth-lthorn,prop:thorns-sym,prop:lodestar-sym}.
\end{proof}

We have the following corollary.

\begin{coro}[{\cite[Conjecture~1.4]{kvec-depth}}] \label{coro:qt-sym-lodestar}
  For any $b \in \mathbb{N}$ and multiset $M$ with elements in $\mathbb{N}$, we have $\tC_{M, b}(q, t) = \tC_{M, b}(t, q)$.
\end{coro}
\begin{proof}
  It is clear that we may write $\tC_{M, b}$ as a sum of $\tC_{a, M' ,b}$ where $a$ runs over elements of $M$ and $M' = M \setminus \{a\}$. We then conclude by \Cref{thm:qt-sym-lodestar}.
\end{proof}

Using our bijections, we also get transparent proofs of some results in \cite{kvec-depth}.

\begin{prop}[{\cite[Proposition~5.3~and~Corollary~5.4]{kvec-depth}}] \label{prop:misc-sym}
  Given a vector $\mathbf{K}$ of natural numbers, we have
  \begin{enumerate}
  \item $\tC_{\mathbf{K} \cdot (a)}(q, t) = \tC_{\mathbf{K} \cdot (b)}(q, t)$ for any $a, b \in \mathbb{N}$; \label{item:res1}
  \item $\tC_{\mathbf{K}}(q, t) = \tC_{\mathbf{K}}(t, q)$ when $\mathbf{K}$ is of the form $(a, b, b, \ldots, b, c)$, hence when $\mathbf{K}$ is of length $3$. \label{item:res2}
  \end{enumerate}
\end{prop}
\begin{proof}
  For \ref{item:res1}, we observe from the proof of \Cref{prop:lodestar-sym} that changing the degree of a lodestar does not change either $\lthorn$ or $\rthorn$. We thus conclude by the fact that the right lodestar of $T$ is always the rightmost internal node in the contour walk order.

  For \ref{item:res2}, we note that it is a natural consequence of \Cref{thm:qt-sym-lodestar}, as there is only one possible sequence $\mathbf{K}$ for $\lukaset_{a, M, c}$ when $M = [b^k]$.
\end{proof}

\begin{rmk}
  We consider the following generating function for the two statistics $\area$ and $\depth$ on \luka{} paths refined by their profile multiset. Let $F(z, q, t; p_0, p_1, \ldots) \equiv F(z, q, t)$ be the generating function defined as
  \[
    F(z, q, t) = \sum_{M \text{ multiset}} z^{|M|} \sum_{P \in \lukaset_{M}} q^{\area(P)} t^{\depth(P)} \prod_{k \geq 0} p_k^{M_k}.
  \]
  Here, $|M|$ denotes the number of elements in the multiset $M$ counted with multiplicity. Through the lens of \Cref{constr:luka-tree}, using \Cref{prop:area-rthorn,prop:depth-lthorn}, by decomposing plane trees at the root and using the standard symbolic method (\emph{cf.} \cite{flajolet}), we have
  \begin{equation}
    \label{eq:fct-eq}
    F(z, q, t) = 1 + \sum_{k \geq 0} p_k \prod_{\ell = 1}^{k} F(z q^{k-\ell} t^{\ell - 1}, q, t).
  \end{equation}
  Here, the first term $1$ stands for the empty tree, and the second term for trees with a root. The degree of the root is given by $k + 1$. For a node in the sub-tree induced by the $\ell$-th child of the root, its $\lthorn$ (resp. $\rthorn$) increments by $\ell - 1$ (resp. $k - \ell$). This explains the substitution of $z$, which counts internal nodes, by $z q^{k - \ell} t^{\ell - 1}$.
\end{rmk}

\bibliographystyle{abbrvnat}
\bibliography{refs}

\end{document}